\documentclass{article}
\usepackage{amssymb}
\usepackage{amsmath}
\usepackage{color}
\usepackage{hyperref}

\def\thebibliography#1{\section*{References}\list
  {[\arabic{enumi}]}{\settowidth\labelwidth{[#1]}\leftmargin\labelwidth
    \advance\leftmargin\labelsep
    \usecounter{enumi}}
    \def\newblock{\hskip .11em plus .33em minus -.07em}
    \sloppy
    \sfcode`\.=1000\relax}

\newcommand{\refbook}[3]{{\sc #1}{\em\ #2}{\ #3}}
\newcommand{\refer}[5]{{\sc #1}{\ #2}{\em\ #3}{\bf\ #4}{\ #5}}

\newtheorem{lem}{Lemma}[section]
\newtheorem{cor}[lem]{Corollary}
\newtheorem{teo}[lem]{Theorem}
\newtheorem{os}[lem]{Remark}
\newtheorem{defi}[lem]{Definition}
\newtheorem{prop}[lem]{Proposition}

\newcommand\sign{\mathop{\rm sign}}
\newcommand{\qed}{\thinspace\null\nobreak\hfill\hbox{\vbox{\kern-.2pt\hrule
 height.2pt depth.2pt\kern-.2pt\kern-.2pt \hbox to2.5mm{\kern-.2pt\vrule
 width.4pt \kern-.2pt\raise2.5mm\vbox to.2pt{}\lower0pt\vtop
 to.2pt{}\hfil\kern-.2pt \vrule
 width.4pt \kern-.2pt}\kern-.2pt\kern-.2pt\hrule height.2pt depth.2pt
 \kern-.2pt}}\par\medbreak}

\newcommand{\R}{\mathbb{R}}

\newcommand{\C}{\mathbb{C}}

\newcommand{\N}{\mathbb{N}}
\newcommand{\Z}{\mathbb{Z}}

\newcommand{\eps}{\varepsilon}

\date{}

\begin{document}


\title{$L^p$ estimates  for Baouendi-Grushin operators}
\author{{G. Metafune \thanks{Dipartimento di Matematica ``Ennio De Giorgi'', Universit\`a del Salento, C.P.193, 73100, Lecce, Italy.
email:  giorgio.metafune@unisalento.it}  \qquad L. Negro  \thanks{Dipartimento di Matematica ``Ennio De Giorgi'', Universit\`a del Salento, C.P.193, 73100, Lecce, Italy.
email:  luigi.negro@unisalento.it}\qquad C. Spina} \thanks{Dipartimento di Matematica ``Ennio De Giorgi'', Universit\`a del Salento, C.P.193, 73100, Lecce, Italy.
email:  chiara.spina@unisalento.it}   }

\maketitle

\begin{abstract}
We prove $L^p$ estimates for the Baouendi-Grushin operator
$\Delta_x+|x|^\alpha \Delta_y$ in $L^p(\R^{N+M})$, $1<p<\infty$, where $x\in \R^N,\ y\in\R^M$.
When $p=2$ more general weights  belonging to the  Reverse H\"older
class $B_2(\R^N)$ are allowed.

\par


\end{abstract}

\section{Introduction}
In this paper we prove $L^p$ estimates for the  Baouendi-Grushin operator
$L=\Delta_x+|x|^\alpha\Delta_y$ in $L^p(\R^{N+M})$, $1<p<\infty$, where $x\in \R^N,\ y\in\R^M$; more specifically, we prove the $L^p$ boundedness
of the operators $ D_{x_ix_j}L^{-1},\
|x|^\alpha D_{y_iy_j}L^{-1}$ and, when $N=1$, also of $|x|^{\frac{\alpha}{2}}D_{xy}L^{-1}$. We use these results to characterize the domain 
 of the operator $L$, denoted by $D_p(L)$, where the solution of the equation $\lambda u-Lu=f$ exists and is unique, for any $f \in L^p(\R^{N+M})$ and $\lambda >0$. In an equivalent way, we describe the domain under which $L$ generates a (analytic and symmetric) semigroup in $L^p(\R^{N+M})$.

\smallskip

When $\alpha$ is an even integer, these $L^p$ estimates are well-known and can be found in the classical paper by Folland, see \cite{folland}. When $\alpha$ is an unresctricted positive real number many results are known on local regularity of the equation $Lu=f$, see for example \cite{FS}, \cite{FGW},\cite{FL} and \cite{garofalo} for unique continuation property. We refere to   \cite{rob-sik} for heat kernel estimates even in a more general context. However, we are not aware of global regularity results for the second derivateves of $u$, with the exception of \cite{wang} where global H\"older regularity is proved and of \cite{kim}, where $L^p$ estimates are proved when $\alpha=1$, $N=1$,  in the half plane $x>0$, for the inhomogenuous problem $Lu=f, u(0,y)=g(y)$. When $g=0$ the estimates in \cite{kim} reduce to ours: even though our results are valid in the whole space, they can be rephrased in the half space $x>0$ when $N=1$ for Dirichlet or Neumann boundary conditions, by considering odd and even (with respect to $x$) functions, respectively.

\smallskip

We  prove $L^p$ estimates through an interpolation theorem in absence of kernels in homogeneous spaces due to Z. Shen, see
\cite[Theorem 3.1]{shen2}, \cite[Theorem 3.14]{auscher1}, and
weighted mean value inequalities for
subsolution of the elliptic equation $Lu=0$ with respect to the balls associated to the subellitic distance defined by the operator, proved in \cite{FS} and \cite{CW}. Some of these results can probably be generalized to the case when $|x|^\alpha$ is replaced by a weight function $\phi(x)$ belonging to the reverse H\"older class $B_p(\R^N)$. This is the case when $p=2$ where the result is not obtained via integration by parts but using maximal results due to \cite{auscher2} for Schr\"odinger operators with $B_2$ potential. However local estimates for subsolutions seem to be known only in special cases and they are crucial in our approach when $p \neq 2$. Another restriction comes from the estimates of the mixed derivatives, that is for the operator $|x|^{\frac{\alpha}{2}}D_{xy}L^{-1}$ where our proof works  when $N=1$ but arbitrary $M$ and relies on a  non standard Rellich type inequality in dimension $1$.

\smallskip

The paper is organized as follows. In Section 2 we define the 
operator in $L^2(\R^{N+M})$ through a form and prove $L^2$ estimates via partial Fourier transform and maximal results on Schr\"odinger operators.  In Section 3 we briefly recall the subelliptic distance associated to $L$ and the main geometrical objects needed in $L^p$ estimates. These last are proved in Section 4, where a separate subsection deals with mixed derivatives.

\smallskip

{\bf Notation} We use $L^p$ for $L^p(\R^{N+M})$, $C_c^\infty$ for
$C_c^\infty (\R^{N+M})$. $L^\infty_c$ stands for the space of all
bounded measurable functions on $\R^{N+M}$ having compact support.
$\mathcal{S}$ is the Schwartz space and $\mathcal S'$  the space of
tempered distributions. 

\section{$L^2$ estimates}\label{Section L^2}
Let $\phi: \R^N\rightarrow [0,+\infty[$ be a nonnegative  continuous function  and set 
$$\Omega_N\subset \R^{N}=\{x\in\R^N:  \phi(x)>0\}, \quad\quad \Omega=\Omega_N\times\R^M.$$
Let $L$ be  the operator defined on smooth functions by 
  $$L=\Delta_x+\phi(x)\Delta_y,$$
 where $x\in \R^N,\ y\in\R^M$. Setting $a=\left(\begin{array}{c|c}
  I_N & 0 \\ 
  \hline
  0 & \phi I_M
 \end{array}\right)=(a_{ij})$ or
\begin{equation} \label{defgamma}
a_{ij}(x,y)=
\begin{cases}
1,&\quad \text{if}\ i=j\leq N,\\
\phi(x),&\quad \text{if}\ N+1\leq i=j\leq N+M,
\end{cases} 
\end{equation}
and $0$ elsewhere, we can write
\begin{equation*}\label{symmetrize}
L={\rm div}(a\,\nabla ) 
\end{equation*}
and, therefore,  $L$ is formally self-adjoint with respect to the Lebesgue measure. 

\begin{os} \label{distribuzione}
Note that $L$ is non-degenerate in the $x$-direction but degenerates in the $y$-direction outside $\Omega$. Accordingly, $\nabla_{x}u$ will denote the distributional gradient (with respect to $x$) of $u$ in the whole space $\R^{N+M}$ and $\nabla_{y}u$ only its distributional gradient (with respect to $y$) in $\Omega$.
\end{os}
We give a formal definition of $L$ through a symmetric form
\begin{defi}\label{def.form.a} 
Consider the sesquilinear form $\mathfrak{a}$ in $L^2$  defined by
\begin{align*}
\mathfrak{a}(u,v)
&:=
\int_{\R^{N+M}}
\Big[\langle \nabla_x u,\nabla_x \overline{v}\rangle+\phi(x)\langle \nabla_y u,\nabla_y \overline{v}\rangle\Big]\,dxdy,\\[1ex]
D(\mathfrak{a})
&:=\{u\in L^2:\ u\in H^1_{loc}(\Omega),\ \nabla_x u,\ \phi^{\frac 1 2}\nabla_y u\in L^2\}.
\end{align*}
\end{defi}
 According to the remark above, we require that the weak gradient $\nabla_y u$ exists only in $\Omega$. 

We summarize in the following lemma the main properties of $a$. Note that, due to the assumptions on $\phi$, $\mathfrak{a}$ is locally uniformly elliptic on $\Omega$.
\begin{lem}\label{Close form Nd}
 $\mathfrak{a}$ is a  nonnegative, symmetric and closed form in $L^2$ and  the following properties hold
\begin{itemize}
%
%
%
%
\item [(i)] if $Q$ is an orthogonal matrix in $\R^M$,  $y_0\in\R^M$ and $I_{Q+y_0} u(x,y)=u(x,Qy+y_0)$, then for every $u,v\in D(\mathfrak{a})$, $I_{Q+y_0}u, I_{Q+y_0} v \in D(\mathfrak{a})$ and
\begin{align*}
\mathfrak{a}(I_{Q+y_0}u, I_{Q+y_0}v)=
\mathfrak{a}(u,v).
\end{align*}
\item [(ii)]if $\phi$ is homogeneous of degree $\alpha \ge 0$, i.e.  $\phi(s x)=s^\alpha\phi(x)$ for $x\in\R^N,\ s>0$, then defining the dilation 
$$I_su(x,y)=u(sx,s^{\frac{2+\alpha}{2}} y)$$ one has for every $u,v\in D(\mathfrak{a})$  $I_su, I_s v \in D(\mathfrak{a})$ and
\begin{align*}
\mathfrak{a}(I_su, I_sv)=s^{2-N-\frac{2+\alpha}{2}M}\,
\mathfrak{a}(u,v);
\end{align*}
\end{itemize}
\end{lem} 
{\sc Proof.}\ 
Clearly, due to the positivity of $\phi$, $\mathfrak{a}$ is a non-negative symmetric form in $L^2$.
The closedness of the form follows easily since $\mathfrak{a}$ is locally uniformly elliptic in $\Omega$ . 
 The proofs of (i) and (ii) follow by a straightforward computation.
\qed
\medspace

Let $-L$ be the operator associated to  $\mathfrak{a}$, that is
\begin{align}\label{Definition operator in L^N}
D(L):=
\left\{u\in D(\mathfrak{a})
\;;\;
\exists v\in L^2\ \text{s.t.}\ \mathfrak{a}(u,w)
=
\int_{\R^{N+M}}v\overline{w}\,d\mu
\quad\forall w\in D(\mathfrak{a})
\right\}, 
\quad 
-Lu:=v.
\end{align}
 The basic properties of $L$ are listed below.
\begin{prop}\label{basic.property}
The operator $-L$ defined in (\ref{Definition operator in L^N}) is nonnegative and selfadjoint. 
Moreover, 
\begin{itemize}
\item [(i)] $C_c^\infty\hookrightarrow D(L)\hookrightarrow  \{u\in L^2\cap W^{2,2}_{\rm loc}(\Omega)\;;\;Lu\in L^2\}$ and for every $u\in C_c^\infty$ 
$$Lu=\Delta_xu+\phi(x) \Delta_yu$$
\item [(ii)] $L$ generates a contractive analytic semigroup $\left\{e^{zL}:\ z\in\C_+\right\}$  in $L^2$.
\item [(iii)] The semigroup $\{e^{tL}:\ t>0\}$ is submarkovian i.e. it is positive and $L^\infty$-contractive. 
\item[(iv)]
If $Q$ is an orthogonal matrix in $\R^M$,  $y_0\in\R^M$ then
$$L=I_{Q+y_0}^{-1}LI_{Q+y_0},\quad  I_{Q+y_0} u(x,y)=u(x,Qy+y_0).$$
\item[(v)] If $\phi$ is homogeneous of degree $\alpha$, then 
\begin{align*}
s^2L=I_s^{-1}LI_s, \quad I_su(x,y)=u(sx,s^\frac{2+\alpha}{2}y), \quad s>0. 
\end{align*}
 
\end{itemize}
\end{prop}
{\sc Proof.}\ (i) is clear by construction and from interior elliptic regularity (see however the proof of Theorem \ref{Teo char D2} for justifying the integration by parts).  The generation property of $L$ follows by  standard results;  the positivity of $e^{tL}$ as well its $L^\infty$-contractivity  is a consequence of  the Beurling-Deny criteria satisfied by the form $\mathfrak{a}$  (see \cite[Corollary 2.18]{ou}). Concerning (iv) and (v), let $u \in D(L)$, $v \in  D(\mathfrak{a})$ and $s>0$. Then 
\begin{align*}
 \mathfrak{a}(I_s u,v)&=s^{2-N-\frac{2+\alpha}{2}M}\mathfrak{a}(u,I_{s^{-1}}v)\\[1ex]
 &=-s^{2-N-\frac{2+\alpha}{2}M} \int_{\R^{N+M}}(Lu){I_{s^{-1}}\overline{v}}\, dxdy=-s^{2}\int_{\R^{N+M}}(I_s Lu)\overline{v} \,dxdy,
\end{align*}

hence $I_su \in D(L)$ and $LI_s u=s^{2}I_s Lu$. Similarly for $I_{Q+y_0}$. \qed

The following Proposition shows that $C_c^\infty$ is dense in $D(L)$ with respect to the graph norm.
\begin{prop}\label{prop density L2}
$C_c^\infty$ is a core for the operator $\left(L,D(L)\right)$ and the form $\mathfrak{a}$.
\end{prop}
{\sc{Proof.}} Since $I-L$ is invertible we have to show that $(I-L)\left(C_c^\infty\right)$ is dense in $L^2$ or, equivalently,  that $(I-L)\left(C_c^\infty\right)^{\perp}=\left\{0\right\}$. To this aim let $v\in L^2$ such that
\begin{align*}
\int_{\R^{N+M}}\left(I-L\right)u\, v\ dx\ dy
=0, \quad \forall u\in C_c^\infty.
\end{align*}
Taking the partial Fourier transform with respect to the $y$ variable and applying Fubini and Plancherel Theorems we get 
\begin{align*}
\int_{\R^{N+M}}\Big[\hat u(x,\xi)-\Delta_x \hat u(x,\xi)+\phi(x)|\xi|^2\hat u(x,\xi)\Big]\, \hat v(x,\xi)\ dx\ d\xi=0, \quad \forall u\in C_c^\infty.
\end{align*}
Choosing $u=A(x)B(y)\in C_c^\infty $ we have $\hat u(x,\xi)=A(x)\hat B(\xi)$ and  
\begin{align}\label{Core2 eq 1}
\int_{\R^{N+M}}\Big[A(x)-\Delta_x A(x)+\phi(x)|\xi|^2A(x)\Big]\,\hat B(\xi) \ \hat v(x,\xi)\ dx\ d\xi=0.
\end{align}
Fix $\xi_0\in\R^M$, $r>0$ and let $w(\xi)=\frac{1}{|B(\xi_0,r)|}\chi_{B(\xi_0,r)}\in L^2(\R^M)$. Let  $(B_n)_n\in C_c^\infty(\R^{M})$ a sequence of test function such that $B_n\to \check{w}$ in $L^2(\R^M)$; then  $\hat B_n\to w$ in $L^2(\R^M)$ and taking  the limit for $n\to\infty$ in Equation \ref{Core2 eq 1} with $\hat B$ replaced by $\hat B_n$ we obtain
\begin{align*}
\frac{1}{|B(\xi_0,r)|}\int_{B(\xi_0,r)}d\xi\int_{\R^{N}}\Big[A(x)-\Delta_x A(x)+\phi(x)|\xi|^2A(x)\Big]\,\hat v(x,\xi)\ dx=0.
\end{align*}

Letting $r\to 0$ and using the Lebesgue Differentiation theorem, we have for a.e. $\xi_0\in\R^M$
\begin{align*}
\int_{\R^{N}}\Big[A(x)-\Delta_x A(x)+\phi(x)|\xi_0|^2A(x)\Big]\,\hat v(x,\xi_0)\ dx=0,
\end{align*}
which by the arbitrariness of $u$ is valid for every $A\in C^\infty_c(\R^N)$. The operator $\Delta_x -\phi(\cdot)|\xi|^2$ is a Schr\"odinger operator in $L^2(\R^N)$ with nonpositive potential $-\phi|\xi|^2$ and  $C^\infty_c(\R^N)$ is dense  in  the domain $D(\Delta_x -\phi(\cdot)|\xi|^2)$  with respect to the graph norm. The last equation  then implies $\hat v(\cdot,\xi_0)=0$  for a.e. $\xi_0\in\R^M$ which proves the required claim.

Since $D(L)$ is dense in  $D\left(L^{\frac 1 2}\right)=D(\mathfrak{a})$, the second statement follows from the first. \qed

In order to prove the main result of this section we recall the definition of $B_p$-weights.
Let $1<p\leq\infty$. Then  $\omega\in B_p(\R^N)$, the class of
the reverse H\"{o}lder weights of order $p$, if $\omega\in
L^p_{loc},\ \omega>0$ a. e. and there exists a positive constant $C$
such the inequality
\begin{equation}  \label{reverse1}
\left(\frac{1}{|B|}\int_B\omega^p\right )^\frac{1}{p}\leq\frac{C}{|B|}\int_B
\omega
\end{equation}
holds  for every  ball $B$. If $p=\infty$, the left
hand side of the inequality above has to be replaced by the
essential supremum of $\omega$ on $B$. The smallest positive
constant $C$ such that (\ref{reverse1}) holds is the $B_p$ constant
of $\omega$. We recall that powers $|x|^\alpha$ belong to $B_\infty (\R^N)$ whenever $\alpha \ge 0$. This is easily seen first considering balls of radius 1 and (and  large centers) and then scaling.

\begin{teo}   \label{L2-estimates}
Assume that $\phi: \R^N\rightarrow [0,+\infty[$ belongs to $B_2(\R^N)$. Then  for every $1\leq i,j\leq N$, $1\leq h,k\leq M$ one has
$$\|D_{x_ix_j}u\|_2+  \|\phi D_{y_hy_k}u\|_2\leq C\|Lu\|_2,\quad u\in D(L).$$
Moreover
$$\|\nabla_xu\|_2+\|\phi^{\frac{1}{2}}\nabla_yu\|_2\leq C(\|Lu\|_2+\|u\|_2),\quad u\in D(L).$$ 
If $\phi\in B_N(\R^N)$ then, for every $1\leq i\leq N$, $1\leq h\leq M$,  one has 
$$\|\phi^\frac{1}{2} D_{x_iy_h}u\|_2\leq C\|Lu\|_2,\quad u\in D(L).$$  
\end{teo}
{\sc Proof.} By Proposition \ref{prop density L2} we may assume that  $u\in C_c^\infty$. Consider the partial Fourier transform with respect to the $y$ variable. Let $v(x,\xi)=\hat{u}(x,\xi)$.  Then, set
$Lu=f$,  we have
$$\Delta_x v(x,\xi)-\phi(x)|\xi|^2v(x,\xi)=\hat{f}(x,\xi)\in L^2.$$
Observe now that, for every fixed $\xi\in\R^M$,  $\Delta_x -\phi(\cdot)|\xi|^2$ is a Schr\"odinger operator in $\R^N$ with potential $\phi|\xi|^2$. Moreover, since $\phi\in B_2(\R^N)$, it immediately follows that $\phi|\xi|^2$ satisfies the reverse H\"older condition with the same constant as $\phi$. By \cite[Theorem 1.1, Corollary 1.3]{auscher2}, we have that
$$|\xi|^4\int_{\R^N}\phi(x)^2|v(x,\xi)|^2\ dx\leq C \int_{\R^N}|\hat{f}(x,\xi)|^2\ dx$$ with a constant $C$ not depending on $\xi$. Integrating the last inequality over $\R^M$, we get
$$\int_{\R^{N+M}}|\xi|^4\phi(x)^2|v(x,\xi)|^2\ dx \ d\xi\leq C \int_{\R^{N+M}}|\hat{f}(x,\xi)|^2\ dx\ d\xi.$$
Since $|\cdot|^2v(x,\cdot)=\widehat{\Delta_y u}(x,\cdot)$ we get, using Fubini theorem and the  Plancherel Theorem  in $\R^M$,
$$\int_{\R^{N+M}}\phi(x)^2|\Delta_y u(x,y)|^2\ dx\ dy\leq C \int_{\R^{N+M}}|f(x,y)|^2\ dx\ dy,$$
which reads as  $\|\phi \Delta_y u\|_2\leq C\|Lu\|_2$; by difference we also get $\|\Delta_x u\|_2\leq C\|Lu\|_2$.

The  Calderon-Zygmund Theorem applied separately to  each variables implies
\begin{align*}
\| D_{x_ix_j} u\|_{L^2\left(\R^N\right)}^2\leq C(N)\| \Delta_x u\|_{L^2\left(\R^N\right)}^2,\quad 1\leq i,j\leq N,\\
\| D_{y_hy_k} u\|_{L^2\left(\R^M\right)}^2\leq C(M)\| \Delta_y u\|_{L^2\left(\R^M\right)}^2,\quad 1\leq h,k\leq M.
\end{align*}
Integrating the previous inequalities (with the last one multiplied by $\phi(x)^2$) respectively  over $\R^M$ and $\R^N$ we get the first claim.

Concerning the gradient estimates, it is enough to observe that, by interpolation, for every $\eps>0$, 
$$\|\nabla_x\|_{L^2(\R^N)}\leq \eps\sum_{i,j=1}^N\left\|D_{x_ix_j}u\right\|_{L^2(\R^N)}+\frac{C}{\eps}\|u\|_{L^2(\R^N)}.$$ 
The estimates for the the first order derivatives with respect to $x$ immediately follow  after  integration over $\R^M$ and by using the first part of the Theorem. 
For the gradient  with respect to $y$, we  start, analogously, from
$$\|\nabla_y\|_{L^2(\R^M)}\leq \eps\sum_{h,k=1}^M\|D_{y_hy_k}u\|_{L^2(\R^M)}+\frac{C}{\eps}\|u\|_{L^2(\R^M)}.$$
Choosing $\eps= \phi(x)^\frac{1}{2}$, the claim follows after the integration over $\R^N$ and by using the first part of the Theorem.

Assume now $\phi\in B_{N}\left(\R^N\right)$ and consider the mixed derivative $D_{x_iy_h}$. Its partial Fourier transform, with respect to the $y$ variable is given by $-i\xi_hD_{x_i}\hat{u}(x,\xi)$. As before $\phi|\xi|^2$ satisfies the $B_N(\R^N)$ reverse  H\"older condition with the same constant as $\phi$. By \cite[Remark Corollary 1.5]{auscher2}, 
$$|\xi|^2\int_{\R^N}\phi(x)|\nabla_x v(x,\xi)|^2\ dx\leq C \int_{\R^N}|\hat{f}(x,\xi)|^2\ dx$$ with a constant $C$ not depending on $\xi$.
Integrating over $\R^M$ and using Plancherel Theorem, we get
$$\int_{\R^{N+M}}\phi(y)|D_{x_iy_h} u(x,y)|^2\ dx\ dy\leq C \int_{\R^{N+M}}|f(x,y)|^2\ dx\ dy.$$
\qed

In the following Theorem we characterize the domain of the operator $L$. We formulate it in the case where $\phi$ belongs  to $ B_2(\R^N) \cap B_N(\R^N)$. If $\phi$ belongs only to $B_2(\R^N)$, the mixed derivatives $\phi^{\frac12} D_{x_iy_h}u$ should not be considered in the statement below.

\begin{teo}\label{Teo char D2} If $\phi \in B_2(\R^N) \cap B_N(\R^N)$ then
the domain of the operator $L$ defined in (\ref{Definition operator in L^N}) satisfies
\begin{align}\label{Char D2}
D(L)=\left\{u\in L^2: \nabla_xu, D_{x_ix_j}u\in L^2, \ \phi^{\frac 1 2}\nabla_y,\  \phi D_{y_hy_k}u,\  \phi^{\frac12} D_{x_iy_h}u \in L^2\right\}.
\end{align}
\end{teo}
{\sc{Proof.}} Let $\tilde D(L)$ be the set defined in the right hand side of equality \eqref{Char D2}. Theorem \ref{L2-estimates} then implies $D(L)\subseteq\tilde D(L)$. To prove the equality it is then enough to prove that the operator $\left(L, \tilde D(L)\right)$ is dissipative since in this case $I-L:\tilde D(L)\to L^2$ is  an injective extension of the resolvent operator $I-L: D(L)\to L^2$ and so both operators must coincide. Let $u\in \tilde D(L)$; then, by the definition, for every compact set  $\omega\subset\subset \Omega$, $u,\Delta u\in L^2\left(\omega\right)$ hence $u\in H^2_{loc}\left(\Omega\right)$. Moreover a section argument (see for example \cite[Theorem 2.1.4]{ziemer}) shows that  for a.e. $x\in\Omega_N$ $u(x,\cdot)\in H^2\left(\R^M\right)$ and 
\begin{align*}
\int_{\R^M}\Delta_yu u\,dy=-\int_{\R^{M}}|\nabla_yu|^2\,dy,\quad \text{for a.e. } x\in\Omega_N.
\end{align*}
Then multiplying by $\phi$, integrating in $x$  and using  Fubini's Theorem we get
\begin{align*}
\int_{\R^{N+M}}\phi(x)\Delta_yu u\,dxdy=-\int_{\R^{N+M}}\phi(x)|\nabla_yu|^2\,dxdy.
\end{align*}
An analogous reasoning applied to the $y$-sections shows that 
\begin{align*}
\int_{\R^{N+M}}\Delta_xu u\,dxdy=-\int_{\R^{N+M}}|\nabla_xu|^2\,dxdy,\quad.
\end{align*}
The last two inequalities imply
\begin{align*}
\int_{\R^{N+M}}Lu\, u\,dxdy=-\int_{\R^{N+M}}\left(|\nabla_yu|^2+\phi(x)|\nabla_yu|^2\right)\,dxdy\leq 0
\end{align*}
which, by the arbitrariness of $u\in \tilde D(L)$, implies the dissipativity of $\left(L,\tilde D(L)\right)$.
\\\qed

The next proposition provides   regularity properties of the solution of the resolvent equation with respect to the $y$ variables.
\begin{prop}  \label{regularity}
Let $u\in D(L)$ be such that $u-Lu=f\in C_c^\infty$. Then for every  multiindex $\alpha$ one has $D^\alpha_y u\in D(L)$ and 
\begin{align*}
D^\alpha_y u-L D^\alpha_y u=D^\alpha_y  f.
\end{align*}
In particular $u$ is smooth in the $y$ variable.
\end{prop} 
{\sc{Proof.}} Let $u\in D(L)$ be such that $u-Lu=f\in C_c^\infty$. Then
\begin{align}\label{Reg eq 2}
\int_{\R^{N+M}} \Big( uv+\langle\nabla_x u,\nabla_x v\rangle+\phi(x) \langle\nabla_y u,\nabla_y v\rangle \Big)\ dxdy=\int_{\R^{N+M}} fv\ dxdy,\quad \text{for every } v\in D(\mathfrak{a}).
\end{align}
For  $h\in\R^{M}$ let $D_h$ be the difference quotient $D_hg(z):=\left(g(x,y+h)-g(x,y)\right) $ and let us take, in the last equation, $v=D_{-h}D_hu\in D(\mathfrak{a})$. Then, since $D_{-h}=D_h^\ast$, one has
\begin{align}\label{Reg eq 1}
\nonumber\int_{\R^{N+M}} \Big( |D_hu|^2+|D_h\nabla_xu|^2 \phi(x)+|D_h\nabla_y u|^2\Big)\ dxdy=\int_{\R^{N+M}} D_hfD_hu\ dxdy\\[1ex]
\leq \|D_h f\|_2\|D_h u\|_2\leq \frac{1}{2}\left(\|D_h f\|_2^2+\|D_h u\|_2^2\right).
\end{align}
In particular for every $\omega\subset\subset\Omega$ there exists some positive constant $C=C(\omega)$ such that
\begin{align*}
\|D_h\nabla u\|_{L^2(\omega)}\leq C|h|\|\nabla f\|_{L^2(\omega)}
\end{align*}
for sufficiently small $h$; this proves that $\nabla u$ is weakly differentiable in $\omega$ in the $y$ variable and that $D_{y_i} u\in H^1_{loc}(\Omega)$. Moreover, if $e_1,\cdots e_M$ is the standard basis of $\R^M$, $t\neq 0$ and $h=t e_i$, then dividing by $t$ both members of equation \eqref{Reg eq 1} and taking the limit for $t\to 0$   we obtain
\begin{align*}
\frac{1}{2}\int_{\R^{N+M}} \Big( |D_{y_i}u|^2+|D_{y_i}\nabla_xu|^2 \phi(x)+|D_{y_i}\nabla_y u|^2\Big)\ dxdy\leq \int_{\R^{N+M}} |D_{y_i} f|^2\ dxdy
\end{align*}
which proves that $D_{y_i}u\in D(\mathfrak{a})$. Let us fix now $v\in C_c^\infty$ ; using  \eqref{Reg eq 2} with $v$ replaced by $D_{-te_i}v$ we get
\begin{align*}
\int_{\R^{N+M}} D_{te_i}fv\ dxdy&=\int_{\R^{N+M}} fD_{-te_i}v\ dxdy\\[1ex]
&=\int_{\R^{N+M}} \Big( uD_{-te_i}v+\langle\nabla_x u,\nabla_x D_{-te_i}v\rangle+\phi(x) \langle\nabla_y u,\nabla_y D_{-te_i}v\rangle \Big)\ dxdy\\[1ex]
&=\int_{\R^{N+M}} \Big( D_{te_i}uv+\langle D_{te_i}\nabla_x u,\nabla_x v\rangle+\phi(x) \langle D_{te_i}\nabla_y u,\nabla_y v\rangle \Big)\ dxdy.
\end{align*}
Dividing by $t$ both members of the last equation and taking the limit for $t\to 0$ we obtain
\begin{align*}
\int_{\R^{N+M}} D_{y_i}fv\ dxdy=\int_{\R^{N+M}} \Big( D_{y_i}uv+\langle D_{y_i}\nabla_x u,\nabla_x v\rangle+\phi(x) \langle D_{y_i}\nabla_y u,\nabla_y v\rangle \Big)\ dxdy.
\end{align*}
 Since  by Proposition \ref{prop density L2} $C_c^\infty$ is  a core for  $\mathfrak{a}$, the arbitrariness of $v$ in the last equation proves that $D_{y_i}u\in D(L)$ and $D_{y_i}u-L(D_{y_i}u)=D_{y_i} f$ which is  the required claim for $|\alpha|=1$. An inductive argument  easily proves the claim  for any multiindex $\alpha$. Moreover, since $D^\alpha_y u=(I-L)^{-1}D^\alpha_y f$, then for some $C=C(\alpha)>0$
\begin{align*}
\|D^\alpha_y u\|_2\leq C \|D^\alpha_y  f\|_2.
\end{align*}
The Sobolev embbeding Theorem then proves that $u$ is $C^\infty$ in the $y$ variable.\\\qed

We end this section by proving a version  of  Kato's inequality adapted to $L$  which will be used for proving $L^p$-estimates.
\begin{prop} \label{Kato's-inequality}
Let $u \in D(L)$ and let us define
\begin{eqnarray*} \sign(u)=
\begin{cases}
0&if\ \ \  u(x)=0\\
u(x)/|u(x)|&if\ \ \  u(x)\neq 0.
\end{cases}
\end{eqnarray*}
Then $|u|$ satisfies the following distributional inequality
$$
-\mathfrak{a}(|u|, \varphi)\geq \int_{\R^{N+M}} \sign(u) Lu\, \phi \, dxdy\quad  \text{for any}\quad  0\leq\varphi\in C_c^\infty.
$$
\end{prop}

{\sc Proof.} 
We suppose first that $u\in C_c^\infty$. If
\begin{equation*}\label{eq: kato's 1}
    u_\epsilon(x)=\sqrt{|u|^2+\epsilon^2}.
\end{equation*}
Then $u_\epsilon\geq |u|$ and 
\begin{equation}\label{eq: kato's 2}
u_\epsilon (a\,\nabla u_\epsilon)=u (a\,\nabla u)
\end{equation}
(here $a$ is the matrix defined in \eqref{defgamma}), then (\ref{eq: kato's 2}) implies that
\begin{align}\label{eq: kato's 3-1}
  |\nabla_x u_\epsilon|&\leq |u||u_\epsilon|^{-1}|\nabla_x u|\leq |\nabla_x u|,\\[1ex]
\nonumber  \phi(x)|\nabla_y u_\epsilon|&\leq |u||u_\epsilon|^{-1}\phi(x)|\nabla_y u|\leq \phi(x)|\nabla_y u|.
\end{align} 
Taking the divergence of (\ref{eq: kato's 2}) we obtain
\begin{equation*}
    u_\epsilon Lu_\epsilon+|\nabla_x u_\epsilon|^2+\phi(x)|\nabla_y u_\epsilon|^2=u Lu+ |\nabla_x u|^2+\phi(x) |\nabla_y u|^2
\end{equation*}
so by (\ref{eq: kato's 3-1})
\begin{equation}\label{eq: kato's 4}
  L u_\epsilon\geq \frac{u}{u_\epsilon} L u.
\end{equation}
Integrating by parts the right hand side of (\ref{eq: kato's 4}),  it follows that
\begin{equation*}\label{eq: kato's 5}
-\mathfrak{a}(u_\epsilon, \varphi)\geq \int_{\R^{N+M}} \frac{u}{u_\epsilon}\, L u\ \varphi\, dxdy,\quad \text{for any}\quad 0\leq\varphi\in C_c^\infty.
\end{equation*}
Letting $\epsilon \to 0$ we get
\begin{equation*}
-\mathfrak{a}(|u|, \varphi)\geq \int_{\R^{N+M}} \sign(u) Lu\  \varphi\,dxdy,\quad \text{for any}\quad 0\leq\varphi\in C_c^\infty.
\end{equation*}
Let now $u\in D(L)$  and let $u_n\in C_c^\infty$ be such that $u_n\to u$ in $D(L)$.   Up to a subsequence, if necessary, we can also suppose that $u_n\to u$ almost everywhere. Since also $u_n\rightarrow u$ in $D(\mathfrak {a})$  by the last inequalities
$$
-\mathfrak{a}(|u_n|, \varphi)\geq \int_{\R^{N+M}} \sign(u_n) Lu_n\ \varphi\,dxdy,\quad \text{for any}\quad 0\leq\varphi\in C_c^\infty,
$$
the claim follows letting $n\rightarrow \infty$. 
\qed

\section{The distance $d$ associated to $L$}\label{section distance} 

Let $\alpha>0$ and let 
 $$L=\Delta_x+|x|^\alpha\Delta_y,$$
 be the self-adjoint operator defined in Section \ref{Section L^2} with $\phi(x)=|x|^\alpha$. In this section  we introduce a natural metric $d$ on $\R^{N+M}$ associated to $L$ and which makes the triple $\left(\R^{N+M}, d, \mathcal{L}\right)$ consisting of $\R^{N+M}$ equipped with the distance $d$ and the Lebesgue measure $\mathcal{L}$, an homogeneous space in the sense of Coifman and Weiss (see \cite{Coif-weiss 1, Coif-weiss 2}).  We refer the reader to \cite{FS, FL} (and references therein) for the proofs of the following results and for further details. 

\begin{defi}
Let $\gamma:[0,T]\to\R^{N+M}$ be an absolutely continuous curve. 
 We say that $\gamma$ is a subunit curve if for a.e. $t\in [0,T]$ one has
\begin{align*}
\langle\dot\gamma(t),\xi\rangle\leq |\xi_x|^2+|x|^\alpha |\xi_y|^2, \quad \text{for every } \xi=(\xi_x,\xi_y)\in\R^{N+M}.
\end{align*}
For every $z_1,z_2\in\R^{N+M}$ we define $d:\R^{N+M}\times\R^{N+M}\to\R^+$ as
\begin{align}\label{distance}
\nonumber d(z_1,z_2)&=\inf\left\{ T\in\R^+: \text{ there exists a subunit curve } \gamma:[0,T]\to\R^{N+M},\, \gamma(0)=z_1, \, \gamma(T)=z_2\right\}\\[1ex]
&=\sup\left\{\psi(z_2)-\psi(z_1):\psi\in W^{1,\infty}\left(\R^{N+M}\right),\;|\nabla_x\psi|^2+|x|^\alpha|\nabla_y\psi|^2\leq 1\right\}.
\end{align}
\end{defi}
We remark that $d$ is a well defined distance and that any couple of point $z_1,z_2\in\R^{N+M}$ can be joined by a subunit curve, see \cite[Section 2, Example 3.6]{FS} and \cite[Definition 2.4]{FL}. A proof of the equality in \eqref{distance} can be found in \cite[Proposition 3.1]{Jerison}.

For $z_0\in\R^{N+M}$, $r>0$ we write  $S(z_0,r):=\{z\in\R^{N+M}:d(z_0,z)<r \}$ to denote the balls of $\R^{N+M}$ with respect to the metric $d$. In the next Proposition we clarify the structure of the metric and define an equivalent system of balls  which are explicit and easier to work with. For $z_0=(x_0,y_0)\in\R^{N+M}$, $r>0$ let us  define the cylindrical set
 \begin{align}\label{Equivalent balls}
 Q(z_0,r):=B(x_0,r)\times B\Big(y_0,r(x_0)\Big),\quad r(x_0):=r(r+|x_0|)^\frac{\alpha}{2}.
 \end{align}

\begin{prop}\label{prop distance}
There exists two positive constants $C_1,C_2>0$ such that the distance function $d$ satisfies for every $z_1=(x_1,y_1),z_1=(x_2,y_2)\in\R^{N+M}$
\begin{align*}
C_1 F(z_1,z_2)\leq d(z_1,z_2)\leq C_2 F(z_1,z_2)
\end{align*}
where 
\begin{align*}
F(z_1,z_2)&=|x_1-x_2|+\left(\frac{|y_1-y_2|}{\left(|x_1|+|x_2|\right)^\frac{\alpha}{2}}\wedge|y_1-y_2|^{\frac{2}{2+\alpha}}\right)
.
\end{align*}
In particular 
\begin{align*}
|S(z_0,r)|\simeq 
\begin{cases}
r^{N+M(1+\frac\alpha 2)},&\quad \text{if}\ r\geq |x_0|,\\[1.5ex]
r^{N+M}|x_0|^{N+M\frac{\alpha}{2}},&\quad \text{if}\ r\leq |x_0|,
\end{cases}
\end{align*}
and the metric balls satisfy the doubling property 
\begin{align*}
|S(z_0,sr)|\leq C s^{N+M(1+\frac\alpha 2)}|S(z_0,r)|, \quad \text{for every } z_0\in\R^{N+M},\, s\geq 1.
\end{align*}
Furthermore there exists a constant  $c>1$ such that for every  $z_0=(x_0,y_0)\in\R^{N+M}$, $r>0$
\begin{align*}
  Q(z_0,c^{-1} r)\subseteq S(z_0,r)\subseteq Q(z_0,cr).
 \end{align*}
 In particular $|S(z_0,r)|\simeq r^N(r+|x_0|)^{M\frac \alpha 2}$ and $\left(\R^{N+M}, d, \mathcal{L}\right)$ is a metric space of homogeneous type.
\end{prop}
{\sc{Proof.}}The first part of the statement is proved in \cite[Proposition 5.1, Corollary 5.2]{rob-sik} (take in \cite{rob-sik} $\delta_1=\delta_1'=0$, $\delta_2=\delta_2'=\alpha/2$, $D=D'=N+M(1+\alpha/2)$). A proof of the second part can be found in  \cite[Proposition 2.7, Example 3.6]{FS} and  \cite[Proposition 1]{FGW}.\qed

\section{$L^p$ estimates} \label{section LP}
Let $1<p<\infty$. In this section we assume that   $\phi(x)=|x|^\alpha$, with $\alpha>0$ and consider therefore the operator   
$$L=\Delta_x+|x|^\alpha\Delta_y$$ in $L^p$ with $x\in \R^N,\ y\in\R^M$. Property (iii) Proposition \ref{basic.property} shows that the symmetric semigroup $(e^{tL})_{t\geq 0}$ generated by $L$ in $L^2$ is submarkovian. Then by standard result (see for example \cite[Chapter 3]{ou})  it induces a consistent family of strongly continuous semigroup on $L^p$  for any $1<p<\infty$, still denoted by $(e^{tL})_{t\geq 0}$. Moreover $(e^{tL})_{t\geq 0}$ extends to a  contractive  holomorphic  semigroup on a sector (see \cite[Theorem 3.13]{ou}).
\begin{defi}
For any $p\in (1,\infty)$ we define the sectorial operator $\left(L, D_p(L)\right)$ as the  generator of the extrapolated semigroup $(e^{tL})_{t\geq 0}$ in $L^p$. We also write $D_2(L)=D(L)$.
\end{defi}
Note that $ D_p(L)\cap D(L)$, being a dense invariant set, is by construction a core for $\left(L, D_p(L)\right)$. 

\medskip
Theorem \ref{L2-estimates} holds in the specific situation since    $|x|^\alpha \in B_\infty(\R^N)$ and we prove that those estimate extend to $1<p<\infty $. 

We recall that $|x|^\alpha$ belongs to $A_t(\R^N)$, the class of Muckenhoupt weights of order $t \ge 1$, whenever $0 \le \alpha < N(t-1)$.  This means that
$$
\left (\frac{1}{|B|}\int_B |x|^\alpha\, dx\right) \left (\frac{1}{|B|}\int_B |x|^{\alpha(1-t')}\, dx \right )^{t-1} \le C
$$
for any ball (or cube) $B$ of $\R^N$, see for example \cite[Chapther 7.3]{duo}. However we need Muckenhoupt weights in $\left(\R^{N+M}, d, \mathcal{L}\right)$  with respect to the metric defined in Section 3. Since $|x|^\alpha$ is independent of $y$ and since the balls $S$ in this space are equivalent to the cylinders $Q(z,r)$ defined in \eqref{Equivalent balls}, which are products of balls in $R^N$ and $\R^M$, respectively, one easily verifies that
$$
\left (\frac{1}{|S|}\int_S |x|^\alpha\, dx  dy\right) \left (\frac{1}{|S|}\int_S |x|^{\alpha(1-t')}\, dx dy \right )^{t-1} \le C
$$
for every ball $S$ (or cylinder)  in  $\left(\R^{N+M}, d, \mathcal{L}\right)$.

A theory on these classes of weights in homogeneous spaces  is presented for example in \cite[Chapter
I]{stromberg} to which we refer for the proofs of the results needed
in what follows. In particular, we recall that Muckenhoupt weights induce doubling measure. The following elementary consequence of the definition is crucial in our approach.
%
\begin{prop}  \label{reverse2}
If $\phi (x,y)=|x|^\alpha$, $t \ge 1$  and  $\alpha <N(t-1)$, there exists $c>0$ such that the inequality
\begin{equation}   \label{At}
\left(\frac{1}{|Q|}\int_Q g\right)^t\leq \frac{c}{\phi (Q)}\int_Q
g^t\phi
\end{equation}
holds for all nonnegative functions $g$ and all  cylinders
$Q$ in $\left(\R^{N+M}, d, \mathcal{L}\right)$ . Here
$$\phi (Q)=\int_Q\phi.$$
\end{prop}
\medskip

The $A_t$ property of  $\phi=|x|^\alpha$, combined with mean value inequalities for Baouendi-Grushin operators  allow us to
characterize the domain of the operator. We prove the
following result.

\begin{teo}   \label{Lp-estimates}
For every $1\leq i,j\leq N$, $1\leq h,k\leq M$ the operators $|x|^\alpha D_{y_hy_k}(I-L)^{-1}$,  $D_{x_ix_j}(I-L)^{-1}$, originally defined in $L^2$,  extend to bounded operators in $L^p$.
\end{teo}

The main tool is the following  result due to Z. Shen, see
\cite[Theorem 3.1]{shen2}, which can be considered as a 
version of the Calder\`{o}n-Zygmund theorem in absence of kernel. The original proof,
where Euclidean balls are used, can be modified to work also for our space $\left(\R^{N+M}, d, \mathcal{L}\right)$. Indeed an improved version of Shen's result in more general homogeneous spaces, which covers
the cases of our interest, can be found in \cite[Theorem 3.14 and
Section V ]{auscher1}.
\begin{teo}  \label{shen}
Let $1\leq p_0<q_0\leq \infty$. Suppose that $T$ is a sublinear
bounded operator on $L^{p_0}$. Suppose moreover that
there exist $\alpha_2>\alpha_1>1$, $c>0$ such that
\begin{equation*}   \label{hp.shen}
\left(\frac{1}{|Q|}\int_Q|Tf|^{q_0}\right)^\frac{1}{q_0}\leq
C\left(\frac{1}{|\alpha_1
Q|}\int_{\alpha_1Q}|Tf|^{p_0}\right)^\frac{1}{p_0}
\end{equation*}
for all  cylinders $Q$ and for all $f\in C^\infty_c$, with
support in $\R^{N+M}\setminus\alpha_2 Q$. Then, for  $p_0\leq
p<q_0$, there exists  a positive constant $C$ such that  for all  $f\in C^\infty_c$
$$\|Tf\|_p\leq C\|f\|_p.$$
\end{teo}

\bigskip

We briefly describe our strategy of proof  of Theorem \ref{Lp-estimates}. We first prove the a-priori estimates for $p\geq 2$ by  applying the above theorem to the  operator
$T=|x|^\alpha D_{y_hy_k}(I-L)^{-1}$ with $p_0=2$,  arbitrary $q_0>2$  and
$\alpha_1=3$, $\alpha_2=4$. Therefore we have to prove that, if  $Q$
is a cylinder and $0\leq f\in C_c^\infty$ has support in
$\R^{N+M}\setminus 4Q$, then $u=(I-L)^{-1}f$ satisfies
$$\left(\frac{1}{|Q|}\int_Q \Big||x|^\alpha D_{y_hy_k} u\Big|^{q_0}\right)^\frac{1}{q_0}\leq C
\left(\frac{1}{|3Q|}\int_{3Q} \Big||x|^\alpha D_{y_hy_k} u\Big|^{2}\right)^\frac{1}{2}$$ for some positive $C$ independent of
$f$. Observe that  $u$ satisfies in $4Q$ the  equation
$$u-L u=\Delta_x u+|x|^\alpha\Delta_y u=0.$$
  Moreover, by Proposition \ref{regularity}, the operator $L$ commutes with the second order derivatives with respect to $y$ and  $v= D_{y_hy_k} u$  satisfies  the same equation in $4Q$. 
  
  To get the a-priori estimates in the case $1<p\leq 2$, we apply Shen's Theorem to the adjoint operator $T^*$. As first step we recall a mean value inequality for subsolution of $L$, that is for functions $v$ which satisfy the inequality $Lv \ge v$ in a weak sense.

\begin{lem}(see \cite[Corollary 5.8]{FS}) \label{prelim1}
There exists a positive constant $C$ such that, if $v$ is a local subsolution of $L$ in $4Q$, 
then
$$ \sup_Q |v|\leq C\left(\frac{1}{|3Q|}\int_{3Q}v^2\right)^\frac{1}{2}.$$
\end{lem}

As in \cite[Theorem 4.1]{CW}, we can deduce the previous mean value inequality also for $0<r<\infty$.  

\begin{lem} \label{prelim1-rgen}
For every $0<r<\infty$, there exists a positive constant $C$ such that, if 
 $v$ is local subsolution of $L$ in $4Q$, then
$$ \sup_Q |v|\leq C\left(\frac{1}{|3Q|}\int_{3Q}|v|^r\right)^\frac{1}{r}.$$
\end{lem}

Now we prove that Lemma \ref{prelim1} holds if we replace the Lebesgue measure with  $|x|^\alpha \, dx$.
\begin{lem}   \label{prelim2}
Fix $0<s<\infty$  and $v$ as in
Lemma \ref{prelim1}. Then 
$$\sup_Q |v|\leq\left(\frac{C}{\phi(3Q)}\int_{3Q}\phi |v|^s\right)^\frac{1}{s}$$
where $C$ depends only on $s,p$ and the $B_p$ constant of $\phi(x)=|x|^\alpha$ and
$$\phi(3Q)=\int_{3Q}\phi.$$
\end{lem}
{\sc{Proof.}} Let $0<s<\infty$ and $Q$ be a cylinder of $\R^{N+M}$. We
fix $t$ as in Proposition \ref{reverse2}. By using Lemma
\ref{prelim1-rgen} with $r=\frac{s}{t}$ and (\ref{At})  we obtain
$$\sup_Q |v|\leq C\left(\frac{1}{|3Q|}\int_{3Q}|v|^\frac{s}{t}\right)^\frac{t}{s}
\leq C\left(\frac{1}{\phi(3Q)}\int_{3Q}\phi|v|^s\right)^\frac{1}{s}.$$
\qed

By combining the estimate in Lemma \ref{prelim2} and the $B_p$
property we deduce the following.
\begin{cor}
\label{quasi-Shen} Let $0<s<\infty$, $1<p<\infty$ and $v$ as
in Lemma \ref{prelim1}. Then 
$$\left(\frac{1}{|Q|}\int_Q(|x|^\alpha|v|^s)^p\right)^\frac{1}{p}\leq\frac{C}{|3Q|}\int_{3Q}|x|^\alpha |v|^s,$$
where $C$ depends only on $s,p$ and the $B_p$ constant of $\phi(x)=|x|^\alpha$.
\end{cor}
{\sc{Proof.}}
Using the $B_p$ property  of $\phi=|x|^\alpha$ and Lemma $\ref{prelim2}$ we
obtain
$$\left(\frac{1}{|Q|}\int_Q(\phi |v|^s)^p\right)^\frac{1}{p}\leq
\left(\frac{1}{|Q|}\int_Q\phi^p\right)^\frac{1}{p}\sup_Q |v|^s\leq C
\left(\frac{1}{|Q|}\int_Q \phi\right)\sup_Q |v|^s\leq
\frac{C}{|3Q|}\int_{3Q}\phi |v|^s.$$
\qed

We can now prove our main result.
\bigskip\\
{\sc{Proof.}}  (Theorem \ref{Lp-estimates}). We first consider the operators $|x|^\alpha D_{y_hy_k}(I-L)^{-1}$.
\smallskip

 Let us preliminary treat the case  $p\geq 2$. Let us fix $q_0>2$ and let $Q$ be a cylinder in $\R^{N+M}$
and $f\in C_c^\infty$   a smooth function with support in $\R^{N+M}\setminus 4Q$. We
set 
\begin{align*}
T=|x|^\alpha D_{y_hy_k}(I-L)^{-1},\quad u=(I-L)^{-1}f, \quad v= D_{y_hy_k}(I-L)^{-1}f.
\end{align*}
 By Theorem \ref{L2-estimates}, $T$ is bounded on $L^2$. Since $f=0$ in $4Q$ and by Proposition \ref{regularity}
$v-Lv=0$ in $4Q$. Combining the last equality with Kato's inequality of  Proposition \ref{Kato's-inequality}, we get
$$-\mathfrak{a}(|v|,\varphi)\geq\int\sign v\ Lv\ \varphi=\int|v|\ \varphi\geq 0\quad\forall\ 0\leq\phi\in C_c^\infty (4Q).$$ It follows that $v$ is a local subsolution of $L$.
By Corollary \ref{quasi-Shen} with $s=2$ and $\alpha$ replaced by $2\alpha$ we have 
$$\left(\frac{1}{|Q|}\int_Q(|x|^{2\alpha}|v|^2)^q\right)^\frac{1}{q}\leq\frac{C}{|3Q|}\int_{3Q}|x|^{2\alpha} |v|^2,\quad 1<q<\infty$$
 or, equivalently,
$$\left(\frac{1}{|Q|}\int_Q(|x|^{\alpha}|v|)^{2q}\right)^\frac{1}{2q}\leq \left(\frac{C}{|3Q|}\int_{3Q}(|x|^{\alpha} |v|)^2\right)^\frac{1}{2},\quad 1<q<\infty.$$
It follows that for  $2\leq q_0<\infty$
\begin{align*}
\left(\frac{1}{|Q|}\int_Q |Tf|^{q_0}\right)^\frac{1}{q_0}&=
\left(\frac{1}{|Q|}\int_Q (|x|^\alpha |D_{y_hy_k}u|)^{q_0}\right)^\frac{1}{q_0}\\[1ex]
&\leq
\left(\frac{C}{|3Q|}\int_{3Q}(|x|^{\alpha} |D_{y_hy_k}u|)^2\right)^\frac{1}{2}=\left(\frac{C}{|3Q|}\int_{3Q}|Tf|^2\right)^\frac{1}{2}.
\end{align*}

By Theorem \ref{shen}, $T$ extends to a  bounded operator in $L^p$ for every $2\leq p<q_0$. Since we can choose $q_0$ arbitrarily,the case
$2\leq p<\infty$ follows.
\smallskip

To treat the case $p< 2$ we
consider  the adjoint operator 
$$T^*=D_{y_hy_k}(I-L)^{-1}|x|^\alpha$$
 which is bounded in $L^2$. By duality, the boundedness of   $T$  in $L^p$ for every $1< p\leq 2$ is equivalent to the boundedness of $T^*$ is bounded in $L^p$ for every $p\geq 2$. As before we fix $q_0>2$ and  we prove that $T^*$ satisfies Shen's assumption for every $2\leq q_0<\infty$.  Let  $Q$ be   a cylinder in $\R^{N+M}$
and $f\in C_c^\infty$ with support in $\R^{N+M}\setminus 4Q$;  set 
$$u=(I-L)^{-1}(|x|^\alpha f), \quad v= D_{y_hy_k}u.$$
 Then $v$ satisfies $Lv=v$ in $4Q$. By arguing as above, $v$ is a local subsolution of $L$ hence Lemma \ref{prelim1}
yields  a positive constant $C$  such that
$$ \sup_Q |v|\leq C\left(\frac{1}{|3Q|}\int_{3Q}v^2\right)^\frac{1}{2}.$$ 
If follows that
\begin{align*}
\left(\frac{1}{|Q|}\int_Q |T^*f|^{q_0}\right)^\frac{1}{q_0}&=
\left(\frac{1}{|Q|}\int_Q  |D_{y_hy_k}u|^{q_0}\right)^\frac{1}{q_0}\\[1ex]
&
\leq  \sup_Q |v|\leq C\left(\frac{1}{|3Q|}\int_{3Q}v^2\right)^\frac{1}{2}=\left(\frac{C}{|3Q|}\int_{3Q} |T^*f|^{2}\right)^\frac{1}{2}.
\end{align*}
and the proof is complete by Theorem \ref{shen} again, applied to $T^*$.
\smallskip

By difference the operator $\Delta_x(I-L)^{-1}$ is bounded on $L^p$ and, integrating with respect to $y$ the classical Calderon-Zygmund estimates in the $x$ variables we deduce the $L^p$ boundedness of $D_{x_ix_j}(I-L)^{-1}$.
\qed

We can now give an explicitly description of the domain $D_p(L)$. 

\begin{teo}\label{teo char D_p}
Let $p\in(1,\infty)$. Then one has 
\begin{align}\label{Char Dp}
D_p(L)=\left\{u\in L^p:\ \nabla_x u, D_{xx}u\in L^p,\ |x|^{\frac{\alpha}{2}}\nabla_y u, |x|^\alpha D_{yy}u\in L^p \right\}.
\end{align}
Moreover  
\begin{align}\label{a-priori estimates p}
\|D_{x_ix_j}u\|_p+ \||x|^\alpha D_{y_hy_k}u\|_p\leq C\|Lu\|_p,\quad u\in  D_p(L).
\end{align}
\end{teo}
{\sc{Proof.}} Let $1<p<\infty$ and let $\tilde D_p(L)$ be the set defined in the right hand side of equality \eqref{Char Dp}.

 Let us preliminary prove that $D_p(L)\subseteq \tilde D_p(L)$. \\ Theorem \ref{Lp-estimates} and the consistency of the resolvent operators in $L^2$ and in $L^p$ imply that  
\begin{equation}\label{non-omo 1}
\|D_{x_ix_j}u\|_p+ \||x|^\alpha D_{y_hy_k}u\|_p\leq C(\|(I-L)u\|_p
\end{equation} 
for any $ u\in \left(I-L\right)^{-1}\left(C_c^\infty\right)$ which is dense in $D_p(L)$ with respect to the graph norm. 
 This implies that \eqref{non-omo 1} extends to $D_p(L)$ proving that $u$ has pure second order distributional derivatives which satisfies $D_{x_ix_j},|x|^\alpha D_{y_hy_k}\in L^{p}$ and that 
 \begin{equation} \label{non-omo}
\|D_{x_ix_j}u\|_p+ \||x|^\alpha D_{y_hy_k}u\|_p\leq C(\|(I-L)u\|_p\leq C\left(\|u\|_p+\|Lu\|_p\right) , \quad  u\in D_p(L).
\end{equation}
As in Theorem \ref{L2-estimates}, an interpolation argument shows  that $\nabla_x u,\ |x|^{\frac{\alpha}{2}}\nabla_y u\in L^p(\R^{N+M})$ i.e. $u\in \tilde D_p(L)$.

To get homogeneous estimates, we use  Proposition \ref{basic.property} (v), and  apply (\ref{non-omo}) to $u(x,y)=v(sx, s^\frac{2+\alpha}{2}y)$, $s>0$ thus obtaining
$$\|D_{x_ix_j}u\|_p+ \||x|^\alpha D_{y_hy_k}u\|_p\leq C(\|Lv\|_p+s^{-2}\|v\|_p).$$
Letting $s$ to infinity  we obtain \eqref{a-priori estimates p}.

To prove that $\tilde D_p(L)=D(L)$, we proceed as  in the proof of Theorem \ref{Char D2} and show that the operator $\left(L, \tilde D_p(L)\right)$ is dissipative.  Let $u\in \tilde D_p(L)$; then the same sectional argument of  Theorem \ref{Char D2}  shows that  for a.e. $x\in\Omega_N$ $u(x,\cdot)\in W^{2,p}\left(\R^M\right)$ and from  \cite{met-spi}
\begin{align*}
\int_{\R^M}\Delta_yu u|u|^{p-2}\,dy=-(p-1)\int_{\R^{M}}|\nabla_yu|^2|u|^{p-2}\,dy,\quad \text{for a.e. } x\in\Omega_N.
\end{align*}
Then multiplying by $\phi$, integrating in $x$  and using  Fubini's Theorem we get
\begin{align*}
\int_{\R^{N+M}}\phi(x)\Delta_yu u|u|^{p-2}\,dxdy=-(p-1)\int_{\R^{N+M}}\phi(x)|\nabla_yu|^2|u|^{p-2}\,dxdy.
\end{align*}
Analogously 
\begin{align*}
\int_{\R^{N+M}}\Delta_xu u|u|^{p-2}\,dxdy=-(p-1)\int_{\R^{N+M}}|\nabla_xu|^2|u|^{p-2}\,dxdy,\quad.
\end{align*}
The last two inequalities imply
\begin{align*}
\int_{\R^{N+M}}Lu \,u|u|^{p-2}\,dxdy=-(p-1)\int_{\R^{N+M}}\left(|\nabla_xu|^2+\phi(x)|\nabla_yu|^2\right)|u|^{p-2}\,dxdy\leq 0
\end{align*}
which, by the arbitrariness of $u\in \tilde D_p(L)$, implies the dissipativity of $\left(L,\tilde D_p(L)\right)$.
\\\qed

 The following proposition shows also that $C_c^\infty$ is a core for  $\left(L,D_p(L)\right)$.
\begin{prop}\label{prop density Lp}
For any $p\in(1,\infty)$, $C_c^\infty$ is a core for the operator $\left(L, D_p(L)\right)$.
\end{prop}
{\sc{Proof.}}
Let $u\in D_p(L)$; we preliminary approximate $u$ with functions in
$D_p(L)$ having compact support in $\R^{N+M}$. Let  $\eta\in C^\infty_c\left(\R^N\right)$ be   a smooth function such that $\chi_{B_1}\leq\eta\leq \chi_{B_2}$
  and, for every $n\in\N$, $x\in\R^N$, define $\eta_n(x)=\eta\left(\frac{x}{n}\right)$. Set $u_n=\eta_n u$. $u_n$ has, by construction, compact support in $x$ and, using  the characterization in \eqref{Char Dp}, one can easily recognize that $u_n\in D_p\left(L\right)$. Lebesgue's Theorem  immediately implies that $u_n,\ |x|^\frac{\alpha}{2}\nabla_y u_n, |x|^\alpha D_{yy}u_n$ tend respectively 
to $u, \ |x|^\frac{\alpha}{2}\nabla_y u, |x|^\alpha D_{yy}u$ in $L^p$. Concerning the $x$-gradient, we have
\begin{align*}
 \|\nabla_x(\eta_n u)-\nabla_x u)\|_p^p\leq&
 \int_{\R^{N+M}}|\eta_n-1|^p|\nabla_x u|^p\,dxdy+\int_{\R^{N+M}}|\nabla_x\eta_n|^p |u|^p\,dxdy\\[1ex]
 \leq&\int_{\R^{N+M}}|\eta_n-1|^p|\nabla u|^p\,dxdy+
 C  n^{-p}\int_{\{n \le |x| \le
 2n\}}|u|^p\,dxdy,
\end{align*}
which tends to 0 by dominated convergence. Similarly $D_{xx} u_n$ tends to $D_{xx} u$ in $L^p$. This proves that $u_n$ tends to $u$ in $D_p(\R^N)$. Using a similar argument with $\eta$ replaced by an analogous  cut-off function $\eta\in C_c^\infty(\R^M)$,  we can approximate $u$ with functions in $D_p(L)$ having compact support also in the $y$-variable.\\
Let us suppose, now, $p<2$ and let $u\in D_p(L)$ (the case $p=2$ is already proved in Proposition \ref{prop density L2}). From the first part of the proof, we can suppose $u$ having compact support. Let $\eta\in C_c^\infty$ such that $\eta=1$ on $\mbox{supp }u$ and, using Proposition \ref{prop density L2}, let $\left(u_n\right)_{n\in \N}$ be a sequence of $C_c^\infty$ functions such that $u_n\to u$ in $D_2(L)$ (and, a fortiori, $\eta u_n\to u$) as $n$ goes to $\infty$. This implies
\begin{align*}
\|L\left(\eta u_n-u\right)\|_p+\|\eta u_n-u\|_p\leq |\mbox{sup }u|^{\frac{1}{p}-\frac 1 2}\Big[\|L\left(\eta u_n-u\right)\|_2+\|\eta u_n-u\|_2\Big]\to 0,\quad \text {as }n\to\infty,
\end{align*}
which proves the required claim for $p<2$.

The proof of the case $p>2$ can be carried out by slightly adapting the arguments used in the proof of Proposition \ref{prop density L2}. We equivalently shows that $(I-L)\left(C_c^\infty\right)$ is dense in $L^p$
 and to this aim let $v\in L^{p'}(\R^{N+M})$ such that
\begin{align*}
\int_{\R^{N+M}}\left(I-L\right)u\, v\ dx\ dy
=0, \quad \forall u\in C_c^\infty.
\end{align*}
Since $1<p'<2$,  the partial Fourier transform of $v(x,\cdot)\in L^{p'}(\R^M)$, with respect to the $y$ variable,  exists as a function in $L^p\left(\R^M\right)$ for a.e. $x\in \R^N$. Therefore, taking in the last equality the  Fourier transform with respect to the $y$ variable and applying Fubini and Plancherel Theorems, we get 
\begin{align*}
\int_{\R^{N+M}}\Big[\hat u(x,\xi)-\Delta_x \hat u(x,\xi)+\phi(x)|\xi|^2\hat u(x,\xi)\Big]\, \hat v(x,\xi)\ dx\ d\xi=0, \quad \forall u\in C_c^\infty.
\end{align*}
Proceeding as in the proof of Proposition \ref{prop density L2} we conclude that $\hat v(\cdot,\xi)=0$  for a.e. $\xi\in\R^M$ which proves the required claim.
\qed

\subsection{Mixed  derivatives}
By using classical covering results, Rellich type inequalities, the previous estimates for the second order derivatives and symmetry arguments, we obtain here  $L^p$ estimates for the mixed second order derivatives when $N=1$. 
To simplify the notation we write $y$ for any of the variables $y_h$, $h=1, \dots , M$.

\begin{teo}   \label{mixed-derivatives}
Let $N=1$. For every  $u\in D_p(L)$
$$\||x|^\frac{\alpha}{2} D_{xy}u\|_p\leq C\|Lu\|_p.$$
\end{teo}

We need a  Rellich type inequality.


\begin{lem} \label{Rellich}   
There exist a positive constant $C$ such that for $u\in C_c^\infty$, such that $u(0,y)=0,\ \ u_x(0, y)=0$ we have
$$\left\|\frac{u}{x^2}\right\|_p\leq C \|Lu\|_p.$$
\end{lem}
{\sc Proof.} 
Let $u\in C_c^\infty$, such that $u(0,y)=0,\ \ u_x(0, y)=0$. Then by  \cite[Theorem 4.2]{Crit Rellich} and \cite[Proposition 3.10]{rellich} 
$$\int_{\R}\left|\frac{u}{x^2}\right|^p\, dx\leq C\int_{\R}|u_{xx}|^p\, dx.$$ Integrating the previous inequality over $\R^M$ and  using Theorem \ref{Lp-estimates},
$$\left\|\frac{u}{x^2}\right\|_p\leq C \|u_{xx}\|_p\leq C(\|Lu\|_p+\||x|^\alpha\Delta_y u\|_p)\leq C \|Lu\|_p.$$
\qed

\begin{os}
The above Rellich inequality uses Theorem \ref{Lp-estimates} to replace the second derivative with respect to $x$ with the operator $L$. However, even its version in dimension 1 (that is for $ D_{xx}$ rather than $L$) is not obvious and probably cannot be obtained by integration by parts (see e.g. \cite{Crit Rellich} where it his shown that Rellich inequalities can be proved for the Laplacian in $ L^p(\R^N)$ when $p<N/2$, a condition which is never verified in dimension 1).
\end{os}

We first prove mixed derivatives estimates assuming that both $u$ and $u_x$ vanish for $x=0$, assuming  Theorem \ref{Lp-estimates}. 
\begin{lem} \label{boundaryAss}If
 If $N=1$ then 
$$\||x|^\frac{\alpha}{2} D_{xy}u\|_p\leq C\|Lu\|_p$$
for every $u\in C_c^\infty$, such that $u(0,y)=0,\  u_x(0, y)=0$.
\end{lem}
{\sc Proof.} 
Let $u\in C_c^\infty$, such that $u(0,y)=0,\ \ u_x(0, y)=0$.  We first prove the required estimate on $\R^+$.
For every $n\in\Z$ let $x_n:=2^n$ and let us consider the following   asymmetric dyadic intervals centred on $x_n$ defined by
\begin{align*}
B(x_n)&:=\left]x_n-\frac{2^n}{4},x_n+2^n\right[=2^{n-2}\,]3,8[,\\[1ex]
A(x_n)&:=\left]x_n-\frac{2^n}{2},x_n+2\,2^n\right[=2^{n-2}\,]2,12[.
\end{align*}
 $\{B(x_n),\,n\in\Z\}$ is a countable covering of $\R^+$ 
such that at most a finite number  $\zeta\in\N$ among the double intervals $\{A(x_n),\,n\in\Z\}$ overlap.

We fix $\vartheta\in C_c^{\infty }(\R)$ such that $0\leq \vartheta \leq 1$, $\vartheta(x)=1$ for $x\in ]3,8[$ and $\vartheta(x)=0$ for $x\notin ]2,12[$. Moreover, for $n\in Z$ we set $\vartheta_n(x)=\vartheta \left(\frac{x}{\rho_n}\right)$, where $\rho_n=\frac{1}{4}|x_n|=2^{n-2}$.

We apply the classical $L^p$  estimates for elliptic operators with constant coefficients
to the function $\vartheta_n u$ and obtain
$$\||x_n|^\frac{\alpha}{2} D_{xy}(\vartheta_n u)\|_p\leq C\|(\vartheta_n u)_{xx}+|x_n|^\alpha\Delta_y(\vartheta_n u)\|_p.$$
By the classical interpolation inequalities for the gradient we get
\begin{align*}
\||x_n|^\frac{\alpha}{2} D_{xy} u\|_{L^p(B(x_n))}&\quad\leq C\Big(\|u_{xx}+|x_n|^\alpha\Delta_y u\|_{L^p(A(x_n))}+\frac{1}{\rho_n}\| u_x\|_{L^p(A(x_n))}\\&\quad+\frac{1}{\rho_n^2}\|u\|_{L^p(A(x_n))}+\frac{|x_n|^\frac{\alpha}{2}}{\rho_n}\|\nabla_y u\|_{L^p(A(x_n))}\Big)\\[1ex]
&\quad\leq C\big(\|\ u_{xx}+|x_n|^\alpha\Delta_y u\|_{L^p(A(x_n))}+\| u_{xx}\|_{L^p(A(x_n))}\\&\quad+\||x_n|^\alpha\Delta_y u\|_{L^p(A(x_n))}+\frac{1}{\rho_n^2}\|u\|_{L^p(A(x_n))}\big).
\end{align*}
Since  
\begin{align*}
\rho_n=\frac{1}{4}|x_n|,\quad\frac{|x_n|}{2}\leq |x|\leq 3|x_n|,\qquad x\in A(x_n),
\end{align*}
 then we get
\begin{align*}
\||x|^\frac{\alpha}{2} D_{xy} u\|_{L^p(B(x_n))}   
&\leq C\Big(\|u_{xx}+|x|^\alpha\Delta_y u\|_{L^p(A(x_n))}+\|u_{xx}\|_{L^p(A(x_n))}\\&\quad+\||x|^\alpha\Delta_y u\|_{L^p(A(x_n))}+\left\|\frac{u}{|x|^2}\right\|_{L^p(A(x_n))}\Big).
\end{align*}
Summing up over $n$, it follows that
\begin{align*}
\||x|^\frac{\alpha}{2} D_{xy} u\|_{L^p\left(\R^+\right)}&\leq C\Big(\|u_{xx}+|x|^\alpha\Delta_y u\|_p+\|u_{xx}\|_p+\||x|^\alpha\Delta_y u\|_p+\left\|\frac{u}{|x|^2}\right\|_p\Big).
\end{align*}
Using Theorem \ref{Lp-estimates} and the Rellich inequality of Lemma \ref{Rellich} we get
\begin{align*}
\||x|^\frac{\alpha}{2} D_{xy} u\|_{L^p\left(\R^+\right)}  
&\leq C\|u_{xx}+|x|^\alpha\Delta_y u\|_p.
\end{align*}
By repeating the same argument taking
\begin{align*}
 x_n=-2^{n},\quad B(x_n)=2^{n-2}\,]-8,-3[,\quad  A(x_n)=2^{n-2}\,]-12,-2[,
 \end{align*}
  we obtain the same estimates in $L^p\left(]-\infty,0[\right)$.
\qed

Next we prove mixed derivatives estimates assuming that either $u$ or $u_x$ vanishes for $x=0$.
\begin{lem}   \label{boundarAssum2}
If  $N=1$ and $p\neq\frac{2}{2-\alpha}$, then
$$\||x|^\frac{\alpha}{2} D_{xy}u\|_p\leq C\|Lu\|_p$$
for every $u\in C_c^\infty$, such that $u(0,y)=0$ or $u_x(0, y)=0$.
\end{lem}
{\sc Proof.} Let $u\in C_c^\infty$, such that $u(0,y)=0$ and let $v(x,y)=\frac{1}{\lambda}u(\lambda  x,y)$. Then $v(0,y)=0$ and $v_x(0,y)=u_x(0,y)$. This implies  that $w=u-v$ satisfies $w(0,x)=w_x(0,x)=0$. Moreover
\begin{align*}
\||x|^\frac{\alpha}{2} D_{xy}v\|_p
&=\lambda ^{-\frac{\alpha}{2}-\frac{1}{p}} \||x|^\frac{\alpha}{2} D_{xy} u\|_p
\end{align*}
and, applying  Theorem \ref{Lp-estimates},
\begin{align*}
\||x|^\frac{\alpha}{2} Lv\|_p
&\leq \|v_{xx}\|_p+\||x|^\alpha\Delta_yv\|_p\\[1ex]
&=\lambda ^{1-\frac{1}{p}}\|u_{xx}\|_p+\lambda ^{-\alpha-1-\frac{1}{p}}  \||x|^\alpha \Delta_y u\|_p\leq C(\lambda)\|Lu\|_p.
\end{align*}
Applying  Lemma \ref{boundaryAss} to $w$ we then have
\begin{align*}
\||x|^\frac{\alpha}{2} D_{xy}u\|_p  &\leq \||x|^\frac{\alpha}{2} D_{xy} w\|_p+\||x|^\frac{\alpha}{2} D_{xy} v\|_p  
\leq C\left(\|Lw\|_p+\||x|^\frac{\alpha}{2} D_{xy} v\|_p\right)\\[1ex]
&\leq C\left(\|Lu\|_p+ \|Lv\|_p+\||x|^\frac{\alpha}{2} D_{xy} v\|_p\right)
\leq C(\lambda)\|Lu\|_p+ C\||x|^\frac{\alpha}{2} D_{xy} v\|_p\\[1ex]
&= C(\lambda)\|Lu\|_p+ C \lambda ^{-\frac{\alpha}{2}-\frac{1}{p}}  \||x|^\frac{\alpha}{2} D_{xy} u\|_p.
\end{align*}
The claim then follows by choosing $\lambda$ large enough such that $C\lambda^{-\frac{\alpha}{2}-\frac{1}{p}}\leq \frac{1}{2}$.

\smallskip

Assume now $u_x(0,y)=0$ and let $v(x,y)=u(\lambda  x,y)$. Then $u(0,y)=v(0,y)$ and $v_x(0,y)=\lambda u_x(0,y)=0$. Moreover
$$||x|^\frac{\alpha}{2} D_{xy} v\|_p=
\lambda ^{1-\frac{\alpha}{2}-\frac{1}{p}} \||x|^\frac{\alpha}{2} D^2_{xy} u\|_p.$$
It follows that $w=u-v$ satisfies $w(0,x)=w_x(0,x)=0$. Hence an analogous argument yields 
\begin{align*}
\||x|^\frac{\alpha}{2} D_{xy}u\|_p\leq C(\lambda)\|Lu\|_p+ C \lambda ^{1-\frac{\alpha}{2}-\frac{1}{p}}  \||x|^\frac{\alpha}{2} D_{xy} u\|_p.
\end{align*}
Choosing $\lambda$ large enough or small enough accordingly to  $1-\frac{\alpha}{2}-\frac{1}{p}>0$ or $1-\frac{\alpha}{2}-\frac{1}{p}<0$ we get the claim for $1-\frac{\alpha}{2}-\frac{1}{p}\neq 0$ or, equivalently, $p\neq\frac{2}{2-\alpha}$.\qed

{\sc Proof.} (Theorem \ref{mixed-derivatives}).Let us suppose, preliminarily, $p\neq\frac{2}{2-\alpha}$ and let  $u\in C_c^\infty$. 
We introduce the operators
$$Pu(x,y)=\frac{u(x,y)+u(-x,y)}{2},\quad\quad  Qu(x,y)=\frac{u(x,y)-u(-x,y)}{2}$$
Observe that 
$$u=Pu+Qu, \quad (Pu)_x(0,y)=(Qu_x)(0,y)=0,$$
 $P$ and $Q$ commute with the second order derivatives and $\|P(Lu)\|_p+\|Q(Lu)\|_p$ is equivalent to $\|Lu\|_p$.  Moreover 
 $$L(Pu)=P(Lu),\quad L(Qu)=Q(Lu).$$
  We can therefore apply the results in Lemma \ref{boundarAssum2} to $Pu$ and $Qu$.
For the mixed second order derivatives we get
\begin{align*}
\||x|^\frac{\alpha}{2} D_{xy} u\|_p&\leq  \||x|^\frac{\alpha}{2} P(D_{xy} u)\|_p+\||x|^\frac{\alpha}{2} Q(D_{xy} u)\|_p = \||x|^\frac{\alpha}{2} D_{xy} (Pu)\|_p+\||x|^\frac{\alpha}{2} Q_{xy}(Qu)\|_p\\[1ex]
&\leq C\left(\|L(Pu)\|_p+\|L(Qu)\|_p\right)= C\left(\|P(Lu)\|_p+ \|Q(Lu)\|_p\right)\leq C\|Lu\|_p.
\end{align*}
By density the proof extends to $u\in D_p(L)$.

Suppose now  $p=\frac{2}{2-\alpha}$. Observe that, by the previous part of the proof, the operator $|x|^\frac{\alpha}{2}D_{xy}(I-L)^{-1}$ is bounded in $L^{p}$ for $p< \frac{2}{2-\alpha}$ and for $p> \frac{2}{2-\alpha}$. The Riesz-Thorin interpolation Theorem then yields the boundedness of $|x|^\frac{\alpha}{2}D_{xy}(I-L)^{-1}$ also for $p= \frac{2}{2-\alpha}$;  the same scaling argument used in the proof of Theorem \ref{teo char D_p} then proves the required claim.
 \qed

\end{document}